\documentclass[11pt]{amsart}

\usepackage{latexsym}
\usepackage{amsfonts}
\usepackage{amssymb}
\usepackage{amsmath}
\usepackage{mathrsfs}
\input xypic
\usepackage{hyperref}

\newcommand{\be}{\begin{equation}}
\newcommand{\ee}{\end{equation}}
\newcommand{\bea}{\begin{eqnarray}}
\newcommand{\eea}{\end{eqnarray}}
\newcommand{\bean}{\begin{eqnarray*}}
\newcommand{\eean}{\end{eqnarray*}}
\newcommand{\brray}{\begin{array}}
\newcommand{\erray}{\end{array}}

\newtheorem{dfn}{Definition}[section]
\newtheorem{thm}[dfn]{Theorem}
\newtheorem{lmma}[dfn]{Lemma}
\newtheorem{ppsn}[dfn]{Proposition}
\newtheorem{crlre}[dfn]{Corollary}
\newtheorem{xmpl}[dfn]{Example}
\newtheorem{rmrk}[dfn]{Remark}

\newcommand{\bdfn}{\begin{dfn}\rm}
\newcommand{\bthm}{\begin{thm}}
\newcommand{\blmma}{\begin{lmma}}
\newcommand{\bppsn}{\begin{ppsn}}
\newcommand{\bcrlre}{\begin{crlre}}
\newcommand{\bxmpl}{\begin{xmpl}}
\newcommand{\brmrk}{\begin{rmrk}\rm}

\newcommand{\edfn}{\end{dfn}}
\newcommand{\ethm}{\end{thm}}
\newcommand{\elmma}{\end{lmma}}
\newcommand{\eppsn}{\end{ppsn}}
\newcommand{\ecrlre}{\end{crlre}}
\newcommand{\exmpl}{\end{xmpl}}
\newcommand{\ermrk}{\end{rmrk}}

\newcommand{\bbz}{\mathbb{Z}}

\author{S. Sundar\\ }
\title{Inverse semigroups and the Cuntz-Li algebras }

\email{sundarsobers@gmail.com}
\address{Currently visiting The Institute of Mathematical Sciences, Chennai.}
\begin{document}
\maketitle 
\begin{abstract}
 In this paper, we apply the theory of inverse semigroups to the  $C^{*}$-algebra $U[\mathbb{Z}]$  considered in \cite{Cuntz}. We show that the $C^{*}$-algebra $U[\mathbb{Z}]$ is generated by an inverse semigroup of partial isometries. We explicity identify  the groupoid $\mathcal{G}_{tight}$ associated to the inverse semigroup and show that $\mathcal{G}_{tight}$ is exactly the same groupoid obtained in \cite{Cuntz-Li}.
\end{abstract}

\section{Introduction}
 Ever since the appearance of the  Cuntz algebras $O_{n}$  and the Cuntz-Krieger algebras $O_{A}$ there has been a great deal of interest in understanding the structure of  $C^{*}$-algebras generated by partial isometries. The theory of graph $C^{*}$-algebras owes much to these examples. It has now been well known  that these algebras admit a groupoid realisation and the groupoid turns out to be r-discrete. Another object that is closely related with an r-discrete groupoid is that of an inverse semigroup. The relationship between r-discrete groupoids and inverse semigroups was already clear from \cite{Renault}. 

 An inverse semigroup $S$ is a semigroup together with an involution $*$ such that for every $s \in S$,  $s^{*}ss^{*}=s^{*}$ . The universal example of an inverse semigroup is the semigroup of partial bijections on a set. Just like one can associate a $C^{*}$-algebra to a group , one can associate a universal $C^{*}$-algebra related with an inverse semigroup $S$ and is denoted $C^{*}(S)$. This universal $C^{*}$-algebra captures the representations of the inverse semigroup ( as partial isometries on a Hilbert space). One can canonically associate an r-discrete groupoid $\mathcal{G}_{S}$ to an inverse semigroup $S$ such that the $C^{*}$-algebra of the groupoid $\mathcal{G}_{S}$ coincides with $C^{*}(S)$. For a more detailed account of inverse semigroups and r-discrete groupoids, we refer to \cite{Alan} and \cite{Ex}. 

Recently, Cuntz and Li in \cite{Cuntz-Li} has introduced a $C^{*}$-algebra associated to every integral domain with only finite quotients. Earlier in \cite{Cuntz}, Cuntz considered the integral domain $\mathbb{Z}$. Let $R$ be an integral domain with only finite quotients. Then the universal algebra $U[R]$ is the universal $C^{*}$-algebra generated by a set of unitaries $\{u^{n}:n \in R\}$ and a set of partial isometries $\{s_{m}:m \in R^{\times}\}$ satisfying certain relations. In \cite{Cuntz-Li}, it was proved that $U[R]$ is simple and purely infinite. Moreover Cuntz and Li obtained a groupoid realisation of it which they later used it to compute the $K$-groups of these algebras for specific integral domains (See \cite{Cuntz-Li-1} and \cite{Cuntz-Li2}).   A concrete realisation of $U[R]$ can be obtained by representing $s_{m}$ and $u^{n}$  on $\ell^{2}(R)$ by 
\begin{align*}
 s_{m}\to S_{m}:\delta_{r}\to \delta_{rm} \\
 u^{n}\to U^{n}:\delta_{r}\to \delta_{r+n}
\end{align*}
Then $U[R]$ is isomorphic to the $C^{*}$-algebra generated by $S_{m}$ and $U^{n}$ ( by the simplicity of $U[R]$). The operator $S_{m}$ is implemented by the multiplication by $m$ (an injection) and $U^{n}$ is implemented by the  addition by $n$( a bijection). Thus it is immediately clear that $U[R]$ is generated by an inverse semigroup of partial isometries. Thus the theory of inverse semigroups should explain some of the results obtained by Cuntz and Li in \cite{Cuntz-Li}. The purpose of this paper is to obtain the groupoid realisation (obtained in \cite{Cuntz-Li}) by using the theory of inverse semigroups. We spell out the details only for the case $R=\mathbb{Z}$ as the analysis for general integral domains with finite quotients is similar. We should also remark that alternate approaches to the Cuntz-Li algebras were considered  in \cite{Ex2} and in \cite{Quigg-Landstad}. We should mention that this paper contains no new results. The point is if one uses the language of inverse semigroups one can obtain a groupoid realisation systematically without having to guess anything about the structure of the Cuntz-Li algebras.
 
Now we indicate the organisation of the paper. In Section $2$, the definition of $U[\mathbb{Z}]$ is recalled and we show that $U[\bbz]$ is generated by an inverse semigroup of partial isometries which we denote it by $T$. In Section $3$, we recall the notion of tight representations of an inverse semigroup, a notion introduced by Exel in \cite{Ex}. We show that the identity representation of $T$ in $U[\mathbb{Z}]$ is in fact tight,  and show that $U[\mathbb{Z}]$ is isomorphic to the $C^{*}$-algebra of the groupoid $\mathcal{G}_{tight}$ (considered in \cite{Ex}) associated to $T$. In Sections $4$ and $5$, we explicity identify the groupoid $\mathcal{G}_{tight}$ which turns out to be exactly the groupoid considered in \cite{Cuntz-Li}. In Section $6$, we show that $U[\mathbb{Z}]$ is simple. In section $7$, we digress a bit to explain the connection between Crisp and Laca's boundary relations and Exel's tight representations of  Nica's inverse semigroup. In the final Section, we give a few remarks of how to adapt the analysis carried out in Sections $1-6$ for a general integral domain. 
 A bit of notation: For non-zero integers $m$ and $n$, we let $[m,n]$ to denote the lcm of $m$ and $n$ and $(m,n)$ to denote the gcd of $m$ and $n$. For a ring $R$, $R^{\times}$ denotes the set of non-zero elements in $R$. 
        
\section{The Regular $C^*$-algebra associated to $\mathbb{Z}$ }

\begin{dfn}[\cite{Cuntz}]
 Let $U[\mathbb{Z}]$ be the universal $C^*$-algebra generated by a set of
unitaries $\{u^{n}:n \in \mathbb{Z}\}$ and a set of isometries $\{s_{m}: m \in
\mathbb{Z}^{\times} \}$ satisfying the following relations.
\begin{equation*}
\begin{split}
 s_{m}s_{n} &= s_{mn} \\
 u^{n}u^{m}&=u^{n+m} \\
 s_{m}u^{n}&=u^{mn}s_{m}\\
\sum_{n \in \mathbb{Z}/(m)}u^{n}e_{m}u^{-n}&=1
\end{split}
 \end{equation*}
where $e_{m}$ denotes the final projection of $s_{m}$.
\end{dfn}
\begin{rmrk}
 Let $\chi$ be a character of the discrete multiplicative group $\mathbb{Q}^{\times}$. Then the universal property of the $C^{*}$-algebra $U[\mathbb{Z}]$ ensures that there exists an automorphism $\alpha_{\chi}$ of the algebra $U[\mathbb{Z}]$ such that $\alpha_{\chi}(s_{m})=\chi(m)s_{m}$ and $\alpha_{\chi}(u^{n})=u^{n}$. This action of the character group of the multiplicative group $\mathbb{Q}^{\times}$ was considered in \cite{Cuntz-Li}.
\end{rmrk}
For $m\neq 0$ and $n \in \mathbb{Z}$, Consider the operators $S_{m}$ and $U^{n}$ defined on $\ell^{2}(\mathbb{Z})$ as follows:
\begin{align*}
 S_{m}(\delta_{r})&=\delta_{rm} \\
 U^{n}(\delta_{r})&=\delta_{r+n}
\end{align*}
Then $s_{m} \to S_{m}$ and $u^{n} \to U^{n}$ gives a representation of the
universal $C^{*}$-algebra $U[\mathbb{Z}]$ called the regular representation 
and its image is denoted by $U_{r}[\mathbb{Z}]$. We begin with a series of Lemmas (highly inspired and adapted from \cite{Cuntz} and from \cite{Cuntz-Li}) which will be helpful in proving that $U[\mathbb{Z}]$ is generated by an inverse semigroup of partial isometries. 

\begin{lmma}
\label{decomposition}
 For every $m,n \neq 0$, one has $e_{m}=\sum_{k \in
\mathbb{Z}/(n)}u^{mk}e_{mn}u^{-mk}$.
\end{lmma}
\textit{Proof.} One has \begin{equation*}
                         \begin{split}
          e_{m}&=s_{m}s_{m}^{*}\\
                &=s_{m}(\sum_{k \in \mathbb{Z}/(n)}u^{k}e_{n}u^{-k})s_{m}^{*}\\
                &=\sum_{k \in
\mathbb{Z}/(n)}s_{m}u^{k}s_{n}s_{n}^{*}u^{-k}s_{m}^{*}\\
                &=\sum_{k \in
\mathbb{Z}/(n)}u^{km}s_{m}s_{n}s_{n}^{*}s_{m}^{*}u^{-km}\\
                &=\sum_{k \in \mathbb{Z}/(n)}u^{km}s_{mn}s_{mn}^{*}u^{-km}\\
                &=\sum_{k \in \mathbb{Z}/(n)}u^{km}e_{mn}u^{-km}.
                \end{split}
                        \end{equation*}
This completes the proof. \hfill $\Box$
\begin{lmma}
\label{product}
 For every $m,n \neq 0$, one has $e_{m}e_{n}=e_{[m,n]}$ where $[m,n]$ denotes
the least common multiple of $m$ and $n$.
\end{lmma}
\textit{Proof.}  Let $c:=[m,n]$ be the lcm of $m$
and $n$.
Then $c=am=bn$ for some $a,b$. Now from Lemma \ref{decomposition}, it follows
that
 \begin{align*}
  e_{m}e_{n}= \sum_{r \in \mathbb{Z}/(a), s \in
\mathbb{Z}/(b)}u^{mr}e_{c}u^{-mr}u^{ns}e_{c}u^{-ns}
 \end{align*}
The product $u^{mr}e_{c}u^{-mr}u^{ns}e_{c}u^{-ns}$ survives if and only if
$mr\equiv ns \mod c$. But the only choice for such an $r$ and an $s$ is when
$r\equiv 0\mod a$ and $s \equiv 0 \mod b$.
[ Reason : Suppose there exists $r$ and $s$ such that $mr\equiv ns \mod c$. Then
$\frac{mr-ns}{c}$ is an integer. That is $\frac{r}{a}-\frac{s}{b}$ is an
integer. Multiplying by $b$, one has that $\frac{br}{a}-s$ and hence
$\frac{br}{a}$ is an integer. But $a$ and $b$ are relatively prime. Hence $a$
divides $r$. Similarly $b$ divides $s$]. Thus $e_{m}e_{n}=e_{c}$. This completes
the proof. \hfill $\Box$

\begin{lmma}
\label{genproduct}
 Suppose $r\neq s$ in $\mathbb{Z}/(d)$ then the projections $u^{r}e_{m}u^{-r}$
and $u^{s}e_{n}u^{-s}$ are orthogonal where $d$ is the gcd of $m$ and $n$.
\end{lmma}
\textit{Proof.} First note that $e_{d}u^{-r}u^{s}e_{d}u^{-s}u^{r}=0$. Hence
$e_{d}u^{-r}u^{s}e_{d}=0$. Now note that
\begin{align*}
 u^{r}e_{m}u^{-r}u^{s}e_{n}u^{-s} & = u^{r}e_{m}e_{d}u^{-r}u^{s}e_{d}e_{n}u^{-s}
~\textrm{ [by Lemma \ref{product}~] }\\
                                  &=
u^{r}e_{m}(e_{d}u^{-r}u^{s}e_{d})e_{n}u^{-s}\\
                                  &=0
\end{align*}
This completes the proof. \hfill $\Box$
\begin{lmma}
\label{genproduct1}
 Let $m,n \neq 0$ be given. Let $d=(m,n)$ and $c=[m,n]$. Suppose $r\equiv s\mod
d$. Let $k$ be such that $k\equiv r\mod m$ and $k \equiv s\mod n$.  Then
$u^{r}e_{m}u^{-r}u^{s}e_{n}u^{-s}=u^{k}e_{c}u^{-k}$.
\end{lmma}
\textit{Proof.} First note that $u^{r}e_{m}u^{-r}=u^{k}e_{m}u^{-k}$ and
$u^{s}e_{n}u^{-s}=u^{k}e_{n}u^{-k}$. The result follows from Lemma
\ref{product}. \hfill $\Box$

\begin{lmma}
\label{another product}
 For $m,n \neq 0$, one has $s_{m}^{*}e_{n}s_{m}=e_{n^{'}}$ where
$n^{'}:=\frac{n}{(n,m)}$.
\end{lmma}
\textit{Proof.} First note that without loss of generality, we can assume that
$m$ and $n$ are relatively prime. Otherwise write $m:=m_{1}d$ and $n:=n_{1}d$
where $d$ is the gcd of $m$ and $n$. Then $(m_{1},n_{1})=1$ and  
\begin{align*}
 s_{m}^{*}e_{n}s_{m}&=s_{m_{1}}^{*}s_{d}^{*}s_{d}s_{n_{1}}s_{n_{1}}^{*}s_{d}^{*}
s_{d}s_{m_{1}}\\
                    &=s_{m_{1}}^{*}e_{n_{1}}s_{m_{1}}
\end{align*}
So now assume $m$ and $n$ are relatively prime. Observe that
$s_{m}^{*}e_{n}s_{m}$ is a selfadjoint  projection.
For
$s_{m}^{*}e_{n}s_{m}s_{m}^{*}e_{n}s_{m}=s_{m}^{*}e_{n}e_{m}s_{m}=s_{m}^{*}e_{m}
e_{n}s_{m}=s_{m}^{*}e_{n}s_{m}$.
Again,
\begin{align*}
 (s_{m}^{*}e_{n}s_{m})^{2}&=s_{m}^{*}e_{n}e_{m}s_{m} \\
                          &=s_{m}^{*}e_{mn}s_{m} ~[\textrm{ by Lemma
\ref{product}} ~]~~ \\
                          &=s_{m}^{*}s_{m}s_{n}s_{n}^{*}s_{m}^{*}s_{m} \\
                          &=e_{n} 
\end{align*}
This completes the proof. \hfill $\Box$
\begin{lmma}
\label{First product with isometry}
 Let $m,n \neq 0$ and $k \in \mathbb{Z}$ be given. If $(m,n)$ does not divide
$k$ then one has $s_{m}^{*}u^{k}e_{n}u^{-k}s_{m}=0$.
 \end{lmma}
\textit{Proof.}
It is enough to show that $x:=e_{n}u^{-k}s_{m}$ vanishes. Thus it is enough to
show that $xx^{*}=e_{n}u^{-k}e_{m}u^{k}e_{n}$. Now Lemma \ref{genproduct}
implies that $xx^{*}=0$. This completes the proof.
\hfill $\Box$

\begin{lmma}
\label{Product with isometry}
 Let $m,n \neq 0$ and $k \in \mathbb{Z}$ be given. Suppose that $d:=(m,n)$
divides $k$. Choose an integer $r$ such that $mr\equiv k \mod n$. Then
$s_{m}^{*}u^{k}e_{n}u^{-k}s_{m} = u^{r}e_{n_{1}}u^{-r}$ where
$n_{1}=\frac{n}{d}$.
\end{lmma}
\textit{Proof.} Now observe that $u^{k}e_{n}u^{-k}=u^{mr}e_{n}u^{-mr}$. Hence
one has
\begin{align*}
 s_{m}^{*}u^{k}e_{n}u^{-k}s_{m}&=s_{m}^{*}u^{mr}e_{n}u^{-mr}s_{m} \\
                               &=u^{r}s_{m}^{*}e_{n}s_{m}u^{-r}\\
                               &=u^{r}e_{n_{1}}u^{-r} ~[\textrm{ by Lemma
\ref{another product}}~]
\end{align*}
This completes the proof. \hfill $\Box$
\begin{rmrk}
\label{semigroup of projections}
 Let $P:=\{u^{n}e_{m}u^{-n}: m \neq 0, n \in \mathbb{Z} \}
\cup \{0\}$. Then the above observations show that $P$ is a commutative
semigroup of projections which is invariant under the map $x\to
s_{m}^{*}xs_{m}$.
\end{rmrk}

The proof of the following proposition is  adapted from
\cite{Cuntz-Li}.

\begin{ppsn}
\label{Inverse semigroup}
 Let $T:=\{s_{m}^{*}u^{n}e_{k}u^{n^{'}}s_{m^{'}}:m,m^{'},k \neq 0, n ,n^{'} \in
\mathbb{Z}\} \cup \{0\}$. Then $T$ is an inverse semigroup of partial
isometries. Let $P:=\{u^{n}e_{m}u^{-n}: m \neq 0, n \in \mathbb{Z} \}
\cup \{0\}$. Then the set of projections in $T$ coincide with $P$. Also the linear span of $T$ is dense in $U[\mathbb{Z}]$.
\end{ppsn}
\textit{Proof.} The fact that $T$ is closed under multiplication follows from
the following calculation.
\begin{align*}
 s_{m}^{*}u^{n}e_{r}u^{-n^{'}}s_{m^{'}}s_{k}^{*}u^{\ell}e_{s}u^{-\ell^{'}}s_{k^{'
}} &=
s_{m}^{*}u^{n}e_{r}u^{-n^{'}}s_{m^{'}}s_{m^{'}}^{*}s_{k}^{*}s_{m^{'}}u^{\ell}e_{
s}u^{-\ell^{'}}s_{k^{'}} \\
  &=
s_{m}^{*}u^{n-n^{'}}u^{n^{'}}e_{r}u^{-n^{'}}e_{m^{'}}s_{k}^{*}s_{m^{'}}u^{\ell}
e_{s}u^{-\ell}u^{\ell-\ell^{'}}s_{k^{'}} \\
   &=
s_{m}^{*}u^{n-n^{'}}\tilde{e}e_{m^{'}}s_{k}^{*}s_{m^{'}}\tilde{f}u^{\ell-\ell^{'
}}s_{k^{'}}~~ [\textrm{~where} ~\tilde{e}=u^{n^{'}}e_{r}u^{-n^{'}} ~\textrm{and~}
\tilde{f}=u^{\ell}e_{s}u^{-\ell}]\\
 &= s_{m}^{*}u^{n-n^{'}}s_{k}^{*}(s_{k}\tilde{e}s_{k}^{*})(s_{k}e_{m^{'}}s_{k}^{*})(s_{m^{'}}\tilde{f}s_{m^{'}}^{*})s_{m^{'}}u^{\ell-\ell^{'}}s_{k^{'}} \\
 &=s_{mk}^{*}u^{kn-kn^{'}}pu^{\ell m^{'}-\ell^{'}m^{'}}s_{k^{'}m^{'}} ~~[ \textrm{~where} ~p:=(s_{k}\tilde{e}s_{k}^{*})(s_{k}e_{m^{'}}s_{k}^{*})(s_{m^{'}}\tilde{f}s_{m^{'}}^{*}) \in P~]
\end{align*}
Thus we have shown that $T$ is closed under multiplication. Clearly $T$ is closed under the involution $*$. Thus the linear span of $T$ is a $*$ algebra containing $s_{m}$ and $u^{n}$ for every $m \neq 0$ and $n \in \mathbb{Z}$. Hence the linear span of $T$ is dense in $U[\mathbb{Z}]$.

 Now we show that every element of $T$ is infact a partial isometry.
Let $v:=s_{m}^{*}u^{n}e_{k}u^{n^{'}}s_{m^{'}}$ be given. Now,
\begin{align*}
 vv^{*}&=s_{m}^{*}u^{n}e_{k}u^{n^{'}}s_{m^{'}}s_{m^{'}}^{*}u^{-n^{'}}e_{k}u^{-n}s_{m} \\
       &=s_{m}^{*}u^{n}(e_{k}u^{n^{'}}e_{m^{'}}u^{-n^{'}}e_{k})u^{-n}s_{m}\\
       &= s_{m}^{*}u^{n}eu^{-n}s_{m} ~[~\textrm{where}~e:=(e_{k}u^{n^{'}}e_{m^{'}}u^{-n^{'}}e_{k}) \in P~]
\end{align*}
Now it follows from Remark \ref{semigroup of projections} that $vv^{*} \in P$. It also shows that the set of projections in $T$ coincides with $P$. This completes the proof. \hfill $\Box$

 The following equality will be used later. Let us isolate it now. 
\begin{equation}
\label{rule of composition}
 s_{m_1}^{*}u^{k_1}s_{n_1}s_{m_2}^{*}u^{k_2}s_{n_2}=s_{m_{1}m_{2}}^{*}u^{m_{2}k_{1}}e_{m_{2}n_{1}}u^{k_{2}n_{1}}s_{n_{1}n_{2}}
\end{equation}

\section{ Tight representations of an  inverse semigroup  }
 Let us recall the notion of tight characters and tight representations from \cite{Ex}.
\begin{dfn}
 Let $S$ be an inverse semigroup with $0$. Denote the set of projections in $S$ by $E$. A character for $E$ is a map $x:E \to \{0,1\}$ such that 
\begin{enumerate}
 \item the map $x$ is a semigroup homomorphism, and
 \item $x(0)=0$.
\end{enumerate}
\end{dfn}
We denote the set of characters of $E$ by $\hat{E}_{0}$. We consider $\hat{E}_{0}$ as a locally compact Hausdorff topological space where the topology on $\hat{E}_{0}$ is the subspace topology induced from the product topology on 
$\{0,1\}^{E}$. 

For a character $x$ of $E$, let $A_{x}:=\{e \in E: x(e)=1\}$. Then $A_{x}$ is a nonempty set satisfying the following properties.
\begin{enumerate}
 \item[(1)]The element $0 \in A_{x}$.
  \item[(2)]If $e \in A_{x}$ and $f \geq e$ then $f \in A_{x}$.
  \item[(3)]If $e,f \in A_{x}$ then $ef \in A_{x}$.
 \end{enumerate}

Any nonempty subset $A$ of $E$ for which $(1),(2)$ and $(3)$ are satisfied is called a filter. Moreover if $A$ is a filter then the indicator function $1_{A}$ is a character. Thus there is a bijective correspondence between the set of characters and filters. A filter is called an ultrafilter if it is maximal. We also call a character $x$ maximal or an ultrafilter if its support $A_{x}$ is maximal.  The set of maximal characters is denoted by $\hat{E}_{\infty}$ and its closure in $\hat{E}_{0}$ is denoted by $\hat{E}_{tight}$.

The following characterization of maximal characters is due to Exel and we refer to \cite{Ex1} for a proof. Let $E$ be an inverse semigroup of projections. Let $e,f \in E$. We say that $f$ intersects $e$ if $fe\neq 0$.
\begin{lmma}
\label{characterization}
 Let $E$ an inverse semigroup of projections with $0$ and $x$ be a character of $E$. Then the following are equivalent.
\begin{enumerate}
 \item The character $x$ is maximal.
 \item The support $A_{x}$ contains every element of $E$ which intersects every element of $A_{x}$.
\end{enumerate}
\end{lmma}
\begin{crlre}
\label{maximality}
 Let $A$ be a unital $C^{*}$-algebra and $E \subset A$ be an inverse semigroup of projections containing $\{0,1\}$. Suppose that $E$ contains a finite set $\{e_{1},e_{2},\cdots,e_{n}\}$ of mutually orthogonal projections such that $\sum_{i=1}^{n}e_{i}=1$. Then for every maximal character $x$ of $E$, there exists a unique $e_{i}$ for which $x(e_{i})=1$.
\end{crlre}
\textit{Proof.} The uniqueness of $e_{i}$ is clear as the projections $e_{1},e_{2},\cdots,e_{n}$ are orthogonal. Now to show the existence of an $e_{i}$ in $A_{x}$, we prove by contradiction. Assume that $e_{i} \notin A_{x}$ for every $i$. Then by Lemma \ref{characterization}, we have that for every $i$, there exists an $f_{i} \in A_{x}$ such that $e_{i}f_{i}=0$. Let $f=\prod f_{i}$. Then $f \in A_{x}$ and thus nonzero and also $fe_{i}=0$ for every $i$. As $\sum_{i}e_{i}=1$, this forces $f=0$. Thus we have a contradiction. \hfill $\Box$

Let us recall the notion of tight representations of semilattices from \cite{Ex} and from \cite{Ex1}. The only semilattice we consider is that of an inverse semigroup of projections or in otherwords the idempotent semilattice of an inverse semigroup. Also our semilattice contains a maixmal element $1$. First let us recall the notion of a cover from \cite{Ex}.
\begin{dfn}
 Let $E$ be an inverse semigroup of projections containing $\{0,1\}$ and $Z$ be a subset of $E$. A subset $F$ of $Z$ is called a cover for $Z$ if given a non-zero element $z \in Z$ there exists an $f \in F$ such that $fz\neq 0$. The set $F$ is called a finite cover if $F$ is finite.
\end{dfn}

The following definition is actually Proposition 11.8 in \cite{Ex}
\begin{dfn}
 Let $E$ be an inverse semigroup of projections containing $\{0,1\}$. A representation $\sigma:E \to \mathcal{B}$ of the semilattice $E$ in a Boolean algebra $\mathcal{B}$ is said to be tight if for every finite cover $Z$ of the interval $[0,x]:=\{z \in E: z \leq x\}$, one has $\sup_{z \in Z} \sigma(z)=\sigma(x)$.
\end{dfn}

Let $A$ be a unital $C^{*}$ algebra and $S$ be an inverse semigroup containing $\{0,1\}$. Let $\sigma:S \to A$ be a unital representation of $S$ as partial isometries in $A$. Let $\sigma(C^{*}(E))$ be the $C^{*}-$subalgebra in $A$ generated by $\sigma(E)$. Then $\sigma(C^{*}(E))$ is a unital, commutative $C^{*}-$algebra and hence the set of projections in it is a Boolean algebra which we denote by $\mathcal{B}_{\sigma(C^{*}(E))}$. We say the representation $\sigma$ is \textbf{tight} if the representation $\sigma:E \to \mathcal{B}_{\sigma(C^{*}(E))}$ is \textbf{tight}.
\begin{lmma}
\label{tightness}
 Let $X$ be a compact metric space and  $E \subset C(X)$ be an inverse semigroup
of projections containing $\{0,1\}$. Suppose that for every finite set of
projections $\{f_{1},f_{2},\cdots,f_{m}\}$ in $E$, there exists a finite set of
mutually orthogonal non-zero projections $\{e_{1},e_{2},\cdots,e_{n}\}$ in $E$ and a 
matrix $(a_{ij})$ such that
\begin{align*}
\sum_{i=1}^{n}e_{i}&=1 \\
f_{i}&=\sum_{j}a_{ij}e_{j}. 
\end{align*}
Then the identity representation of $E$ in $C(X)$ is tight.
\end{lmma}
\textit{Proof.} Let $e \in E\{0\}$ be given and let $F$ be a finite cover for the interval $[0.e]$. Without loss of generality, we can assume that $e=1$ (Just cut everything down by $e$). Let $F:=\{f_{1},f_{2},\cdots,f_{m}\}$. Then by the hypothesis there exists a finite set of mutually orthogonal projections $\{e_{1},e_{2},\cdots,e_{n}\}$ and a matrix $(a_{ij})$ such that $f_{i}=\sum_{j}a_{ij}e_{j}$ and $\sum_{i}e_{i}=1$. For a given $j$, let $A_{j}:=\{i:a_{ij} \neq 0\}$. Since $F$ covers $C(X)$, it follows that for every $j$, $A_{j}$ is nonempty. In otherwords, given $j$, there exists an $i$ such that $f_{i} \geq e_{j}$. Thus $f:=\sup_{i}f_{i} \geq e_{j}$ for every $j$. Hence $f \geq \sup_{j}e_{j} =1$. This completes the proof. \hfill $\Box$

In the next proposition, $T$ denotes the inverse semigroup associated to $U[\mathbb{Z}]$ in Proposition \ref{Inverse semigroup}.
\begin{ppsn}
\label{tightness of the identity representation}
 The identity representation of $T$ in $U[\mathbb{Z}]$ is tight. 
\end{ppsn}
\textit{Proof.} We apply Lemma \ref{tightness}. Let $\{u^{r_{1}}e_{m_{1}}u^{-r_{1}},u^{r_{2}}e_{m_{2}}u^{-r_{2}},\cdots,u^{r_{k}}e_{m_{k}}u^{-r_{k}}\}$ be a finite set of non-zero projections in $P$. By Lemma \ref{decomposition}, it follows that each $f_{i}:=u^{r_{i}}e_{m_{i}}u^{-r_{i}}$ is a linear combination of $\{u^{s}e_{c}u^{-s}: s \in \mathbb{Z}/(c)\}$ where $c$ is the lcm of $m_{1},m_{2},\cdots,m_{k}$. Then Lemma \ref{tightness} implies that the identity representation of $T$ in $U[\mathbb{Z}]$ is tight. This completes the proof. \hfill $\Box$.

Now we will show that the $C^{*}-$algebra of the groupoid $\mathcal{G}_{tight}$ of the inverse semigroup $T$ is isomorphic to the algebra $U[\mathbb{Z}]$. First let us recall the construction of the groupoid $\mathcal{G}_{tight}$ considered in \cite{Ex}. Let $S$ be an inverse semigroup with $0$ and let $E$ denote its set of projections. Note that $S$ acts on $\hat{E}_{0}$ partially. 
For $x \in \hat{E}_{0}$ and  $s \in S$, define $(x.s)(e)=x(ses^*)$. Then
\begin{itemize}
 \item The map $x.s$ is a semigroup homomorphism, and 
 \item $(x.s)(0)=0$.
\end{itemize}
But $x.s$ is nonzeo if and only if $x(ss^*)=1$. For $s \in S$, define the domain
and range of $s$ as 
\begin{align*}
 D_s:&=\{x \in \hat{E}_{0}: x(ss^*)=1\} \\
 R_s:&=\{x \in \hat{E}_{0}: x(s^*s)=1\} 
\end{align*}
Note that both $D_s$ and $R_s$ are compact and open. Moreover $s$ defines a
homoemorphism from $D_s$ to $R_s$ with $s^*$ as its inverse. Also observe that
$\hat{E}_{tight}$ is invariant under the action of $S$.

Consider the transformation groupoid $\Sigma:=\{(x,s):x\in D_s\}$ with the
composition and the inversion being given by:
 \begin{align*}
  (x,s)(y,t):&=(x,st) \textrm{~if~} y=x.s\\
   (x,s)^{-1}:&=(x.s,s^*) 
 \end{align*}
Define an equivalence relation $\sim$ on $\Sigma$ as $(x,s)\sim (y,t)$ if $x=y$
and if there exists an $e \in E$ such that $x \in D_e$ for which $es=et$. Let
$\mathcal{G}=\Sigma/\sim$. Then $\mathcal{G}$ is a groupoid as the product and
the inversion respects the equivalence relation $\sim$. Now we describe a
toplogy on $\mathcal{G}$ which makes $\mathcal{G}$ into a topological groupoid.

For $s \in S$ and $U$ an open subset of $D_s$, let $\theta(s,U):=\{[x,s]: x \in
U \}$. We refer to \cite{Ex} for the proof of the following two propositions. We
denote $\theta(s,D_s)$ by $\theta_s$. Then $\theta_s$ is homeomorphic to $D_s$
and hence is compact, open and Hausdorff.

\begin{ppsn}
 The collection $\{\theta(s,U): s \in S, U \textrm{~open in ~} D_s \}$ forms a
basis for a topology on $\mathcal{G}$. The groupoid $\mathcal{G}$ with this
topology is a topological groupoid whose unit space can be identified with
$\hat{E}_{0}$.
Also one has the following.
\begin{enumerate}
 \item For $s,t \in S$, $\theta_s \theta_t=\theta_{st}$,
 \item For $s \in S$, $\theta_{s}^{-1}=\theta_{s^{*}}$, and
 \item The set $\{1_{\theta_s}: s \in T\}$ generates the $C^*$ algebra
$C^*(\mathcal{G})$.
\end{enumerate}
\end{ppsn}

We define the groupoid $\mathcal{G}_{tight}$ to be the reduction of the groupoid
$\mathcal{G}$ to $\hat{E}_{tight}$. In \cite{Ex}, it is shown that the representation $s \to 1_{\theta_{s}} \in C^{*}(\mathcal{G}_{tight})$ is tight and any tight representation factors through this universal one.

\begin{ppsn}
\label{main proposition}
 Let $T$ be the inverse semigroup associated to $U[\mathbb{Z}]$ in Proposition \ref{Inverse semigroup}. Let $\mathcal{G}_{tight}$ be the tight groupoid associated to $T$. Then $U[\mathbb{Z}]$ is isomorphic to $C^{*}(\mathcal{G}_{tight})$.
\end{ppsn}
\textit{Proof.} Let $t_{m},v^{n}$ denote the images of $s_{m},u^{n}$ in $C^{*}(\mathcal{G}_{tight})$. The universality of the $C^{*}-$algebra $C^{*}(\mathcal{G}_{tight})$ together with Proposition \ref{tightness of the identity representation} implies that there exists a homomorphism $\rho:C^{*}(\mathcal{G}_{tight}) \to U[\mathbb{Z}]$ such that $\rho(t_{m})=s_{m}$ and $\rho(v^{n})=u^{n}$. 

Note that the mutually orthogonal set of projections $\{u^{r}e_{m}u^{-r}: r\in \mathbb{Z}/(m)\}$ cover $T$. Since the representation of $T$ in $C^{*}(\mathcal{G}_{tight})$ is tight, it follows that $\sum_{r}v^{r}t_{m}t_{m}^{*}v^{-r}=1$. Now the universal property of $U[\mathbb{Z}]$ implies that there exists a homomorphism $\sigma:U[\mathbb{Z}]\to C^{*}(\mathcal{G}_{tight})$ such that $\sigma(s_{m})=t_{m}$ and $\sigma(u^{n})=v^{n}$. Now it is clear that $\rho$ and $\sigma$ are inverses of each other. This completes the proof. \hfill $\Box$

In the next two sections, we identify the groupoid $\mathcal{G}_{tight}$ explicitly.

\section{Tight characters of the inverse semigroup $T$}
In this section, we determine the tight characters of the inverse semigroup $T$ defined in Proposition \ref{Inverse semigroup}. Let us recall a few ring theoretical notions. We denote the set of strictly positive integers by $\mathbb{N}^{+}$. Consider the directed set $(\mathbb{N}^{+},\leq)$ where we say $m \leq n $  if $m|n$. If $m|n$ then there exists a natural map from $\mathbb{Z}/(n)$ to $\mathbb{Z}/(m)$. The inverse limit of this system is called the profinite completion of $\mathbb{Z}$ and is denoted $\hat{\mathbb{Z}}$.  In
other words,
\begin{displaymath}
 \hat{\mathbb{Z}}:=\{(r_{m})\in \prod_{m\in \mathbb{N}^{+}}\mathbb{Z}/(m): r_{mk}\cong r_{m}~~ mod ~m \}
\end{displaymath}
 Also $\hat{\mathbb{Z}}$ is a compact ring with the subspace topology induced by the product topology on $\prod \mathbb{Z}/(m)$. Also $\mathbb{Z}$ embedds naturally in $\hat{\mathbb{Z}}$. We also need the easily verifiable fact that the kernel of the $m^{th}$ projection $r=(r_{m})\to r_{m}$ is in fact $m\hat{\mathbb{Z}}$. 

For $r \in
\hat{\mathbb{Z}}$, define a character $\xi_{r}:P \to \{0,1\}$ by the following
formula:
\begin{align*}
\xi_{r}(u^{n}e_{m}u^{-n}):&=\delta_{r_{m},n} \\
\xi_{r}(0):&=0
\end{align*}
In the above formula, the Dirac-delta function is over the set $\mathbb{Z}/(m)$. Thus $\delta_{r_{m},n}=1$ if and only if $r_{m}\equiv n \mod m$.

\begin{ppsn}
\label{tight characters}
 The map $r \to \xi_{r}$ is a topological isomorphism from $\hat{\mathbb{Z}}$ to
$\hat{P}_{tight}$
\end{ppsn}
\textit{Proof.} First let us check that for $r \in \hat{\mathbb{Z}}.$, $\xi_{r}$ is in fact a character and is maximal. Consider an element $r \in \hat{\mathbb{Z}}$. Let $e:=u^{n_{1}}e_{m_{1}}u^{-n_{1}}$ and $f:=u^{n_{2}}e_{m_{2}}u^{-n_{1}}$ be given. Let $d:=(m_{1},m_{2})$ and $c:=[m_{1},m_{2}]$. Suppose $\xi_{r}(e)=\xi_{r}(f)=1$. Then $r_{m_{1}}\equiv n_{1} \mod m_{1}$ and $r_{m_{2}} \equiv n_{2} \mod m_{2}$. Moreover, $r_{c}\equiv r_{m_{i}} \mod m_{i}$ for $i=1,2$. Thus $ef=u^{r_{c}}e_{c}u^{-r_{c}}$ by Lemma \ref{genproduct1} Hence by definition $\xi_{r}(ef)=1$. Now suppose $\xi_{r}(e)=1$ and $e \leq f$. Then by Lemma \ref{genproduct} and Lemma \ref{genproduct1}, it follows that $m_{2}$ divides $m_{1}$ and $r_{m_{2}}\equiv r_{m_{1}} \equiv n_{1} \equiv n_{2} \mod m_{2}$. Hence $\xi_{r}(f)=1$. By definition $0$ is not in the support of $\xi_{r}$. Thus we have shown that the support of $\xi_{r}$ is a filter or in other words $\xi_{r}$ is a character.

 Now we claim $\xi_{r}$ is maximal. This follows from the observation that for every $m \in \mathbb{N}^{+}$, the set of projections $\{u^{n}e_{m}u^{-n}: n \in \mathbb{Z}/(m)\}$  are mutually orthogonal. Thus if $\xi$ is a character then for every $m$ there exists at most one $r_{m}$ for which $\xi(u^{r_{m}}e_{m}u^{-r_{m}})=1$. This implies that if $\xi$ is a character which contains the support of $\xi_{r}$ then $\xi=\xi_{r}$.

Now let $\xi$ be a maximal character of $P$. Then by Corollary \ref{maximality} and by the observation in the previous paragraph, it follows that for every $m$ there exists a unique $r_{m}$ such that $\xi(u^{r_{m}}e_{m}u^{-r_{m}})=1$. Now let $k$ be given. Since both $u^{r_{m}}e_{m}u^{-r_{m}}$ and $u^{r_{mk}}e_{mk}u^{-r_{mk}}$ belong to the support of $\xi$, it follows that the product $u^{r_{m}}e_{m}u^{-r_{m}}u^{r_{mk}}e_{mk}u^{-r_{mk}}$ does not vanish. Then by Lemma \ref{genproduct}, it follows that $r_{mk}\equiv r_{m} \mod m$. Thus $r=(r_{m}) \in \hat{\mathbb{Z}}$ and the support of $\xi_{r}$ is contained in the support of $\xi$. Thus again by the observation in the preceeding paragraph, it follows that $\xi=\xi_{r}$.

It is clear from the definition that the map $r \to \xi_{r}$ is one-one and continuous. As $\hat{\mathbb{Z}}$ is compact, it follows that the range of the map $r \to \xi_{r}$ which is $\hat{P}_{\infty}$ is also compact. Hence $\hat{P}_{\infty}=\hat{P}_{tight}$. Thus we have shown that $r\to \xi_{r}$ is a one-one and onto continuous map from $\hat{\mathbb{Z}}$ to $\hat{P}_{tight}$. Since $\hat{\mathbb{Z}}$ is compact, it follows that the above map is in fact a homeomorphism. This completes the proof. \hfill $\Box$

From now on we will simply write $r(e)$ in place of $\xi_{r}(e)$ if $r\in \hat{\mathbb{Z}}$ and $e\in P$.

\section{The groupoid $\mathcal{G}_{tight}$ of the inverse semigroup $T$}
Let us recall a few ring theoretical constructions. Consider the directed set $(\mathbb{N}^{+},\leq)$ where the partial order $\leq$ is defined by $m \leq n$ if $m$ divides $n$. For $m \in \mathbb{N}^{+}$, let $\mathcal{R}_{m}:=\hat{\mathbb{Z}}$. Let $\phi_{m\ell,m}:\mathcal{R}_{m} \to \mathcal{R}_{\ell m}$ be the map defined by mulitplication by $\ell$. Then $\phi_{m\ell,m}$ is only an additive homomorphism and it does not preserve the multiplication. We let $\mathcal{R}$ be the inductive limit of $(\mathcal{R}_{m},\phi_{m\ell,m})$. Then $\mathcal{R}$ is an abelian group and $\hat{\mathbb{Z}}$ is a subgroup of $\mathcal{R}$ via the inclusion $\mathcal{R}_{1}\subset \mathcal{R}$. Note that $\mathcal{R}$ is a locally compact Hausdorff space. Moreover  the group $P_{\mathbb{Q}}:=\big\{\begin{bmatrix}
                          1 & 0 \\
                          b & a                                                                                                                                                                                                                                                                                                                                                                                                                                                                                                                                                                                                                                                                                                                                                                                   
                                                                                                                                                                                                                                                                                                                                                                                                                                                                                                                                                                                                                                                                                                                                                                                             \end{bmatrix}: a \in \mathbb{Q}^{\times},b \in \mathbb{Q} \big \}$
 acts on $\mathcal{R}$ by affine transformations. The action is descibed explicitly by the following formula. For $x \in \mathcal{R}_{p}$
\begin{displaymath}
\begin{split}
 \begin{bmatrix}
  1 & 0 \\
  \frac{n}{m^{'}} & \frac{m}{m^{'}}
 \end{bmatrix}x = mx+np \in \mathcal{R}_{m^{'}p}
\end{split}
\end{displaymath}
One can check that the above formula defines an action of $P_{\mathbb{Q}}$ on $\mathcal{R}$. We need the following lemma which has already appeared in \cite{Ex2}. We recall the proof for completeness.
 \begin{lmma}
 Let $a:=\frac{n}{m^{'}}$ and $b:=\frac{m}{m^{'}}$. Then
$s_{m^{'}}^{*}u^{n}s_{m}$ depends only on $a$ and $b$. 
\end{lmma}
\textit{Proof.} Suppose $\frac{n_{1}}{m_{1}^{'}}=\frac{n_{2}}{m_{2}^{'}}$ and $\frac{m_{1}}{m_{1}^{'}}=\frac{m_{2}}{m_{2}^{'}}$. Then $n_{1}m_{2}^{'}=n_{2}m_{1}^{'}$ and $m_{1}m_{2}^{'}=m_{1}^{'}m_{2}$. Now, we have
\begin{equation*}
\begin{split}
 s_{m_{1}^{'}}^{*}u^{n_{1}}s_{m_1}&=s_{m_{1}^{'}}^{*}s_{m_{2}}^{*}s_{m_{2}}u^{n_{1}}s_{m_1} \\
                                  &=s_{m_{2}^{'}}^{*}s_{m_{1}}^{*}s_{m_{1}^{'}}^{*}s_{m_{1}^{'}}s_{m_{2}}u^{n_{1}}s_{m_1} \\
                                  &=s_{m_{2}^{'}}^{*}s_{m_{1}}^{*}s_{m_{1}^{'}}^{*}u^{n_{1}m_{2}m_{1}^{'}}s_{m_{1}^{'}}s_{m_1}s_{m_2}\\
                                  &=s_{m_{2}^{'}}^{*}s_{m_{1}^{'}}^{*}s_{m_{1}}^{*}u^{n_{1}m_{2}^{'}m_{1}}s_{m_{1}}s_{m_{1}^{'}}s_{m_2}\\
                                  &=s_{m_{2}^{'}}^{*}s_{m_{1}^{'}}^{*}u^{n_{1}m_{2}^{'}}s_{m_{1}^{'}}s_{m_2}\\
                                  &=s_{m_{2}^{'}}^{*}s_{m_{1}^{'}}^{*}u^{n_{2}m_{1}^{'}}s_{m_{1}^{'}}s_{m_2}\\
                                  &=s_{m_{2}^{'}}^{*}u^{n_{2}}s_{m_{1}^{'}}^{*}s_{m_{1}^{'}}s_{m_2}\\
                                  &=s_{m_{2}^{'}}^{*}u^{n_{2}}s_{m_{2}}
\end{split}
\end{equation*}
This completes the proof. \hfill $\Box$

Now we explicitly identify the groupoid $\mathcal{G}_{tight}$ associated to the inverse semigroup $T$. When we consider transformation groupoids, we  consider only right actions. Thus we let $P_{\mathbb{Q}}$ act on the right on the space $\mathcal{R}$ by defining $x.g=g^{-1}x$ for $x \in \mathcal{R}$ and $g \in P_{\mathbb{Q}}$. We show that that groupoid $\mathcal{G}_{tight}$ of the inverse semigroup $T$ is isomorphic to the restriction of the transformation groupoid $\mathcal{R} \times P_{\mathbb{Q}}$ to the closed subset $\hat{\mathbb{Z}}$. (Here we consider $P_{\mathbb{Q}}$ as a discrete group.) Let us begin with a lemma which will be useful in the proof.

\begin{lmma}
 \label{Omitting projections}
In $\mathcal{G}_{tight}$ one has $[(r,s_{{m}^{'}}^{*}u^{n^{'}}e_{k}u^{n}s_{m})]=[(r,s_{{m}^{'}}^{*}u^{n+n^{'}}s_{m})]$
\end{lmma}
\textit{Proof.} First observe that $[(r,s_{m^{'}}^{*})][(r.s_{m^{'}}^{*},u^{n^{'}}e_{k}u^{n}s_{m})]=[(r,s_{m^{'}}^{*}u^{n^{'}}e_{k}u^{n}s_{m})]$. Thus it is enough to consider the case $m^{'}=1$. Now let $s:=u^{n^{'}}e_{k}u^{n}s_{m}$ , $t:=u^{n+n^{'}}s_{m}$ and $e:=u^{n^{'}}e_{k}u^{-n^{'}}$.. Now observe that $ss^{*}:=ett^{*}$. Hence if $r(ss^{*})=1$ then $r(tt^{*})=1$ and $r(e)=1$. Moreover $es=et$. Thus $[(r,s)]=[(r,t)]$. This completes the proof. \hfill $\Box$.

\begin{thm}
\label{main theorem}
 Let $\phi:\mathcal{R}\times P_{\mathbb{Q}}|_{\hat{\mathbb{Z}}} \to \mathcal{G}_{tight}$ be the map defined by 
\begin{displaymath}
\phi\Big( \big(r,\begin{bmatrix}
                                                                           1 & 0 \\
                                                                         \frac{k}{m} & \frac{n}{m}
                                                                          \end{bmatrix}
\big)\Big)=
[(r,s_{m}^{*}u^{k}s_{n})]
\end{displaymath}
 Then $\phi$ is a topological groupoid isomorphism.
\end{thm}
\textit{Proof.}

\underline{The map $\phi$ is well defined.} 

Let $(r,\begin{bmatrix}
                                                                           1 & 0 \\
                                                                         \frac{k}{m} & \frac{n}{m}
                                                                          \end{bmatrix}
\big)$ be an element in $\mathcal{R} \times P_{\mathbb{Q}}|_{\hat{\mathbb{Z}}}$ . Then we have $mr-k=ns$ for some $s \in \hat{\mathbb{Z}}$. Now we need to show that $r(s_{m}^{*}u^{k}e_{n}u^{-k}s_{m})=1$. By Lemma \ref{Product with isometry}, it follows that $s_{m}^{*}u^{k}e_{n}u^{-k}s_{m}=u^{r_{n}}e_{n_{1}}u^{-r_{n}}$ where $n_{1}:=\frac{n}{(n,m)}$. Thus 
\begin{align*}
                   r(s_{m}^{*}u^{k}e_{n}u^{-k}s_{m})&=r(u^{r_{n}}e_{n_{1}}u^{-r_{n}})\\
                                                    &=\delta_{r_{n_1},r_n} \\
                                                    &= 1 ~[\textrm{ Since ~} r_{n}=r_{n_1} ~\textrm{in~} \mathbb{Z}/(n_{1})]                                                                                                                                                                                                                                                                                                                                            
                                                                                                                                                                                                                                                                                                                                                              \end{align*}

\underline{ Surjectivity of $\phi$:}

First let us show that if $[(r,s_{m}^{*}u^{k}s_{n})] \in \mathcal{G}_{tight}$ then $\Big(r,\begin{bmatrix}
                                                                           1 & 0 \\
                                                                         \frac{k}{m} & \frac{n}{m}
                                                                          \end{bmatrix} \Big) \in \mathcal{R}\times P_{\mathbb{Q}}|_{\hat{\mathbb{Z}}}$. Consider an element $[(r,v:=s_{m}^{*}u^{k}s_{n})]$ in $\mathcal{G}_{tight}$.                                              
 Then $r(vv^{*})=1$ and $vv^{*}:=s_{m}^{*}u^{k}e_{n}u^{-k}s_{m}$. Now Lemma \ref{First product with isometry} and \ref{Product with isometry} implies that $(m,n)|k$.  Let $s$ be an integer such that $ms\equiv k\mod n$. Again Lemma \ref{Product with isometry} implies that $vv^{*}=u^{s}e_{n_{1}}u^{-s}$ where $n_{1}:=\frac{n}{(n,m)}$. Now $r(vv^{*})=1$ implies that $r_{n_{1}}\equiv s\mod n_{1}$. But $r_{n}\equiv r_{n_1}\mod n_{1}$ ( as $r\in \hat{\mathbb{Z}}$). Thus we have $r_{n}\equiv s\mod n_{1}$. This in turn implies that $mr_{n}\equiv ms\equiv k\mod n$. Hence $mr-k \in n\hat{\mathbb{Z}}$. Hence $\Big(r,\begin{bmatrix}
                                                                           1 & 0 \\
                                                                         \frac{k}{m} & \frac{n}{m}
                                                                          \end{bmatrix} \Big) \in \mathcal{R}\times P_{\mathbb{Q}}|_{\hat{\mathbb{Z}}}$. Now the surjectivity of $\phi$ follows from Lemma \ref{Omitting projections}.

\underline{Injectivity of $\phi$:}

 Now suppose $[(r,s_{m_{1}}^{*}u^{k_{1}}s_{n_{1}})]=[(r,s_{m_2}^{*}u^{k_2}s_{n_2})]$. Then by definition there exists a projection of the form $e:=u^{r_{p}}e_{p}u^{-r_{p}}$ such that $e(s_{m_1}^{*}u^{k_1}s_{n_1})=e(s_{m_2}^{*}u^{k_2}s_{n_2})$. Consider a character $\chi$ of the discrete group $\mathbb{Q}^{*}$. Let $\alpha_{\chi}$ be the automorphism of the algebra $U[\mathbb{Z}]$ such that $\alpha_{\chi}(u^{n}))=u^{n}$ and $\alpha_{\chi}(s_{m})=\chi(m)s_{m}$. 
\begin{align*}
 \chi(\frac{n_{1}}{m_1})e(s_{m_1}^{*}u^{k_1}s_{n_1})&=\alpha_{\chi}\big(e(s_{m_1}^{*}u^{k_1}s_{n_1})\big) \\
                                                    &=\alpha_{\chi}\big(e(s_{m_2}^{*}u^{k_2}s_{n_2})\big)\\
                                                    &=\chi(\frac{n_2}{m_2})e(s_{m_2}^{*}u^{k_2}s_{n_2})\\
                                                    &=\chi(\frac{n_2}{m_2})e(s_{m_1}^{*}u^{k_1}s_{n_1})
                                    \end{align*}
Since $e(s_{m_1}^{*}u^{k_1}s_{n_1})\neq 0$, it follows that $\chi(\frac{n_1}{m_1})=\chi(\frac{n_2}{m_2})$ for every character $\chi$ of the discrete, multiplicative group $\mathbb{Q}^{*}$. Thus $\frac{n_1}{m_1}=\frac{n_2}{m_2}$.

Since $e(s_{m_1}^{*}u^{k_1}s_{n_1})=e(s_{m_2}^{*}u^{k_2}s_{n_2})$ in $U_{r}[\mathbb{Z}]$ and $\frac{n_1}{m_1}=\frac{n_2}{m_2}$, it follows immediately that $\frac{k_1}{m_1}=\frac{k_2}{m_2}$. Thus we have shown that $\phi$ is injective.

\underline{The map $\phi$ is a homeomorphism.}

First we show $\phi$ is continuous. Let $(r_{n},g_{n})$ be a sequence in $\mathcal{R}\times P_{\mathbb{Q}}|_{\hat{\mathbb{Z}}}$ converging to $(r,g)$. Since we are considering $P_{\mathbb{Q}}$ as a discrete group, we can without loss of generality assume that $g_{n}=g$ for every $n$. Then, from Lemma \ref{tight characters}, it follows that $\phi(r_{n},g_{n})$ converges to $\phi(r,g)$.

For an open subset $U$ of $\hat{\mathbb{Z}}$ and $g:=\begin{bmatrix}
                                                      1 & 0 \\
                                                      \frac{k}{m} & \frac{n}{m}
                                                     \end{bmatrix}$, consider the open set 
\begin{displaymath}
 \theta(U,g):= \{(r,g): r \in U ~\textrm{and~}r.g \in \hat{\mathbb{Z}}\}.
\end{displaymath}
 Then the collection $\{\theta(U,g):U \overbrace{\subset}^{open} \hat{\mathbb{Z}},~g \in P_{\mathbb{Q}}\}$ forms a basis for $\mathcal{R}\times P_{\mathbb{Q}}|_{\hat{\mathbb{Z}}}$. Moreover $\phi(\theta(U,g))=\theta(U,s_{m}^{*}u^{k}s_{n})$. 
Hence $\phi$ is an open map. Thus we have shown that $\phi$ is a homeomorphism.

\underline{$\phi$ is a groupoid morphism.}

First we show that $\phi$ preserves the source and range. By definition $\phi$ preserves the range. Let $\Big(r,g:=\begin{bmatrix}
                                                                           1 & 0 \\
                                                                         \frac{k}{m} & \frac{n}{m}
                                                                          \end{bmatrix} \Big) \in \mathcal{R}\times P_{\mathbb{Q}}|_{\hat{\mathbb{Z}}}$
be given. Let  $v:=s_{m}^{*}u^{n}s_{n}$. Since $r.g \in \hat{\mathbb{Z}}$, it follows that there exists $t \in \hat{\mathbb{Z}}$ such that $mr-k=nt$. We need to show that $\xi_{r}.v=\xi_{t}$. (Just to keep things clear we write $\xi_{r}$ for the character determined by $r$). It is enough to show that the support of $\xi_{t}$ and that of $\xi_{r}.v$ coincide. But then both the characters are maximal and thus it is enough to show that the support of $\xi_{t}$ is contained in the support of $\xi_{r}.v$. Thus, suppose that $\xi_{t}(u^{\ell}e_{s}u^{-\ell})=1$. Then $t_{ns}\equiv t_{s}\equiv \ell \mod s$. This implies $mr_{ns}-k\equiv nt_{ns} \equiv n\ell \mod ns$. Thus $mr_{ns}\equiv k+n\ell \mod ns$. Let $n_{1}:=\frac{ns}{(ns,m)}$. Now observe that
\begin{align*}
 (\xi_{r}.v)(u^{\ell}e_{s}u^{-\ell})&=\xi_{r}(vu^{\ell}e_{s}u^{-\ell}v^{*})\\
                                    &=\xi_{r}(s_{m}^{*}u^{k}s_{n}u^{\ell}e_{s}u^{-\ell}s_{n}^{*}u^{-k}s_{m})\\
                                   &=\xi_{r}(s_{m}^{*}u^{k+n\ell}e_{ns}u^{-(k+n\ell)}s_{m})\\
                                    &=\xi_{r}(u^{r_{ns}}e_{n_1}u^{-r_{ns}}) ~[\textrm{~By Lemma \ref{Product with isometry}}~]\\
                                    &= \delta_{r_{ns},r_{n_1}}\\
                                   &=1 ~[\textrm{~Since~} r_{ns}=r_{n_1} ~in~\mathbb{Z}/(n_1)] 
\end{align*}
Thus we have shown that the support of $\xi_{t}$ is contained in the support of $\xi_{r}.v$ which in turn implies that $\xi_{t}=\xi_{r}.v$. Hence $\phi$ preserves the source.

Now we show that $\phi$ preserves multiplication. Let $\gamma_{i}:=(r_{i},\begin{bmatrix}
                                                                           1 & 0 \\
                                                                         \frac{k_i}{m_i} & \frac{n_i}{m_i}
                                                                          \end{bmatrix})$ for $i=1,2$.
Since $\phi$ preserves the range and source, it follows that if $\gamma_{1}$ and $\gamma_{2}$ are composable, so do $\phi(\gamma_{1})$ and $\phi(\gamma_{2})$. Observe that
\begin{align*}
 \phi(\gamma_{1})\phi(\gamma_{2})&=[(r_{1},s_{m_1}^{*}u^{k_1}s_{n_1}s_{m_2}^{*}u^{k_2}s_{n_2}]\\
                                 &=[r_{1},s_{m_{1}m_{2}}^{*}u^{m_{2}k_{1}}e_{m_{2}n_{1}}u^{k_{2}n_{1}}s_{n_{1}n_{2}}] ~\big(\textrm{ Eq. \ref{rule of composition}}~\big)\\
                                 &=[r_{1},s_{m_{1}m_{2}}^{*}u^{m_{2}k_{1}+n_{1}k_{2}}s_{n_{1}n_{2}}]~\big(\textrm{Lemma \ref{Omitting projections}}~\big)\\
                                 &=\phi(\gamma_{1}\gamma_{2})
\end{align*}

It is easily verifiable  that $\phi$ preserves inversion. This completes the proof. \hfill $\Box$

\begin{rmrk}
Combining Proposition \ref{main proposition} and Theorem \ref{main theorem}, we obtain that $U[\mathbb{Z}]$ is isomorphic to $C^{*}(\mathcal{R}\times P_{\mathbb{Q}}|_{\hat{\mathbb{Z}}})$ which is Remark 2 in page 17 of
 \cite{Cuntz-Li}.
\end{rmrk}
\section{Simplicity of $U[\mathbb{Z}]$}
 First we  recall a few definitions from \cite{Ren2009}. Let $\mathcal{G}$ be an r-discrete, Hausdorff and  locally compact topological groupoid. Let $\mathcal{G}^{0}$ be its unit space. We denote the source and range maps by $s$ and $r$ respectively. The arrows of $\mathcal{G}$ define an equivalence relation on $\mathcal{G}^{0}$ as follows:
\begin{displaymath}
 x \sim y \textrm{~if there exists~} \gamma \in \mathcal{G} \textrm{~such that~} s(\gamma)=x \textrm{~and~} r(\gamma)=y
\end{displaymath}
A subset $E$ of $\mathcal{G}^{0}$ is said to be invariant if   the orbit of $x$ is contained in $E$ whenever $x \in E$. For $x\in \mathcal{G}^{0}$, define the isotropy group at $x$ denoted $\mathcal{G}(x)$ by $\mathcal{G}(x):=\{\gamma \in \mathcal{G}:~s(\gamma)=r(\gamma)=x\}$.

A groupoid $\mathcal{G}$ is said to be
\begin{itemize}
 \item topologically principal if the set of $x \in \mathcal{G}^{0}$ for which $\mathcal{G}(x)=\{x\}$ is dense in $\mathcal{G}^{0}$.
  \item minimal if the only non-empty open invariant subset of $\mathcal{G}^{0}$ is $\mathcal{G}^{0}$.
\end{itemize}

We need the following theorem. We refer to \cite{Ren2009} for a proof.
\begin{thm}
\label{simplicity}
 Let $\mathcal{G}$ be an r-discrete, Hausdorff and locally compact topological groupoid. If $\mathcal{G}$ is topologically principal and minimal then $C_{red}^{*}(\mathcal{G})$ is simple.
\end{thm}

\begin{ppsn}
The $C^{*}$-algebra $U[\mathbb{Z}]$ is simple.
\end{ppsn}
\textit{Proof.} Let $\mathcal{G}$ denote the groupoid $\mathcal{R}\times P_{\mathbb{Q}}|_{\hat{\mathbb{Z}}}$. Since the group $P_{\mathbb{Q}}$ is solvable, it is amenable and thus by Proposition 2.15 of \cite{Muhly}, it follows that the full groupoid $C^{*}$-algebra $C^{*}(\mathcal{G})$ is isomorphic to the reduced algebra $C^{*}_{red}(\mathcal{G})$. Now we apply Theorem \ref{simplicity} to complete the proof. 

First let us show $\mathcal{G}$ is minimal.  
               Let $U$ be a non-empty open invariant subset of $\mathcal{G}^{0}$. For $m=(m_{1},m_{2},\cdots,m_{n}) \in (\mathbb{Z}\backslash \{0\})^{n}$ and $k \in \mathbb{Z}$, let 
\begin{displaymath}
 U_{m,k}:=\{r\in \hat{\mathbb{Z}}:~r_{m_{i}}\equiv k \mod m_{i}\}
\end{displaymath}
Then, by the Chinese remainder theorem, it follows that the collection $\{U_{m,k}\}$ (where $m$ varies over $(\mathbb{Z} \backslash \{0\})^{n}$ (we let $n$ vary too) and $k\in \mathbb{Z}$) is a basis for the topology on $\hat{\mathbb{Z}}$. Also observe that for a given $m$, $\bigcup_{k \in \mathbb{Z}} U_{m,k}=\hat{\mathbb{Z}}$. Moreover the translation matrix $\begin{bmatrix}
                                   1 & 0\\
                                   k_{1}-k_{2} & 1                                                                                                                                                                                                                                                                                                                                                                                 
                                                  \end{bmatrix}$ 
maps $U_{m,k_{1}}$ onto $U_{m,k_{2}}$. Now since $U$ is non-empty and open, there exists an $m$ and a $k_{0}$ such that $U_{m,k_{0}} \subset U$. But since $U$ is invariant, it follows that $U_{m,k} \subset U$ for every $k \in \mathbb{Z}$. Thus $\bigcup_{k\in \mathbb{Z}}U_{m,k}\subset U$. This forces $U=\hat{\mathbb{Z}}$. This completes the proof. \hfill $\Box$

Now we show $\mathcal{G}$ is topologically principal. Let \begin{displaymath}
E:=\{ r \in \hat{\mathbb{Z}}: r \neq 0 ,r_{p^{i}}=0 ~\forall i,  \textrm{~except for finitely many primes }~p\}                                                           
                                                          \end{displaymath}
If one identifies $\hat{\mathbb{Z}}$ with $\displaystyle \prod_{p ~prime}\hat{\mathbb{Z}_{p}}$ then it is clear that $E$ is dense in $\hat{\mathbb{Z}}$. Now let $r \in E$ be given. We claim that $\mathcal{G}(r)=\{r\}$. Suppose $r.\begin{bmatrix}
1 & 0  \\                                                                                                                                                                                                                       
\frac{k}{m} & \frac{n}{m}                                                                                                                                                                                                                         \end{bmatrix}=r$. Then $mr-k=nr$. But  $r_{p}=0$ except for finitely many primes. Thus it follows that $k$ is divisible by infinitely many primes which forces $k=0$. Now $mr=nr$ and $r\neq 0$ implies $m=n$. Thus $\mathcal{G}(r)=\{r\}$. This proves that $\mathcal{G}$ is topologically principal. This completes the proof. \hfill $\Box$

\section{Nica-covariance, tightness and boundary relations}
In this section, we digress a bit to understand some of the results in \cite{Nica1},\cite{Crisp} and in \cite{Laca1} from the point of view of  inverse semigroups. 
Let us recall the notion of quasi-lattice ordered groups considered by Nica in \cite{Nica1}. Let $G$ be a discrete group and $P$ a subsemigroup of $G$ containing the identity $e$. Also assume that $P\cap P^{-1}=\{e\}$. Then $P$ induces a left-invariant partial order $\leq$ on $G$ defined by $x \leq y$ if and only if $x^{-1}y \in P$. The pair $(G,P)$ is said to be quasi-lattice ordered if the following conditions are satisfied.
\begin{enumerate}
 \item[(1)] Any $x \in PP^{-1}$ has a least upper bound in $P$, and
 \item[(2)] If $s,t \in P$ have a common upper bound in $P$ then $s,t$ have a least common upper bound.
\end{enumerate}
If $s,t \in P$ have a common upper bound in $P$ then we denote the least upper bound in $P$ by $\sigma(s,t)$. It is easy to show that $s,t \in P$ have a common upper bound if and only if $s^{-1}t \in PP^{-1}$. Let us recall the Wiener-Hopf representation from \cite{Nica1}. Consider the representation $W:P \to B(\ell^{2}(P))$ defined by 
\begin{equation*}
 W(p)(\delta_{a}):=\delta_{pa}
\end{equation*}
where $\{\delta_{a}:a \in P\}$ denotes the canonical orthonormal basis of $\ell^{2}(P)$. Note that for $s \in P$, $W(s)$ is an isometry and $W(s)W(t)=W(st)$ for $s,t \in P$. For $s \in P$, let $M(s)=W(s)W(s)^{*}$ then
\begin{eqnarray}
 M(s)M(t) &=& \left\{\begin{array}{ll}
                        M(\sigma(s,t)) &\mbox{if $s$ and $t$ have common upper bound in $P$}\\
                        0 & \mbox{otherwise}.
                        \end{array} \right.
\end{eqnarray}
Let $\mathcal{N}:=\{W(s)W(t)^{*}:s,t \in P\}\cup \{0\}$. Then Equation $(5)$ of Proposition 3.2 in \cite{Nica1} implies that $\mathcal{N}$ is an inverse semigroup of partial isometries. The following definition is due to Nica.
\begin{dfn}[\cite{Nica1}]
 Let $(G,P)$ be a quasi-lattice ordered group. An isometric representation $V:P \to B(\mathcal{H})$ on a Hilbert space $\mathcal{H}$ (i.e. $V(t)^{*}V(t)=1$ for $t \in P$, $V(e)=1$ and $V(s)V(t)=V(st)$ for every $s,t \in P$) is said to be Nica-covariant if  the following holds
\begin{eqnarray}
 L(s)L(t) &=& \left\{\begin{array}{ll}
                        L(\sigma(s,t)) &\mbox{if $s$ and $t$ have common upper bound in $P$}\\
                        0 & \mbox{otherwise}.
                        \end{array} \right.
\end{eqnarray}
where  we set $L(t)=V(t)V(t)^{*}$.
In other words a Nica-covariant representation of $(G,P)$ is nothing but a unital representation of the inverse semigroup $\mathcal{N}$ which sends $0$ to $0$.
\end{dfn}
 Let us say a Nica-covariant representation is tight if the corresponding representation on $\mathcal{N}$ is tight. Now one might ask what are the tight representations of the inverse semigroup $\mathcal{N}$? We prove that tight representations are nothing but Nica-covariant representations satisfying the boundary relations considered by Laca and Crisp in \cite{Crisp}. This fact  is implicit in \cite{Crisp} and it is in fact explicit if one applies Theorem 13.2 of \cite{Ex1}. The author believes that it is worth recording this connection and we do this in the next proposition. 

First let us fix a few notations. A finite subset $F$ of $P$ is said to cover $P$ if given $x \in P$ there exists $y \in F$ such that $x$ and $y$ have a common upper in $P$. Let 
\begin{equation*}
 \mathcal{F}:=\{F \subset P: \text{$F$ is finite and covers $P$}\}
\end{equation*}

\begin{ppsn}
 Let $(G,P)$ be a quasi-lattice ordered group. Consider a  Nica-covariant representation $V:P \to B(\mathcal{H})$. Then $V$ is tight if and only if for every $F \in \mathcal{F}$, one has $\prod_{t\in F}(1-V(t)V(t)^{*})=0$.
\end{ppsn}
\textit{Proof.} Consider a Nica-covariant representation $V:P \to B(\mathcal{H})$. Suppose that $V$ is tight. Let $F \in \mathcal{F}$ be given. Note that $F$ covers $P$ if and only if $\{M(t):t \in F\}$ covers the set of projections in $\mathcal{N}$. Now the tightness of $V$ implies that $\displaystyle \sup_{t \in F}V(t)V(t)^{*}=1$. This is equivalent to saying that $\prod_{x \in F}(1-V(t)V(t)^{*})=0$. Thus we have the implication '$\Rightarrow$'.

Let $V$ be a Nica-covariant representation for which $\prod_{t \in F}(1-V(t)V(t)^{*})=0$ for every $F \in \mathcal{F}$. We denote the set of projections in $\mathcal{N}$ by $E$. Then $E:=\{M(t): t\in P\} \cup \{0\}$. Let $\{M(t_{1}),M(t_{2}),\cdots,M(t_{n})\} \subset [0,M(t)]$ be a finite cover. Then $M(t_{i}) \leq M(t)$ for every $i$. But this is equivalent to the fact that $t \leq t_{i}$. 

We claim that $ \{t^{-1}t_{i}:i=1,2,\cdots,n\}$ covers $P$. Let $s \in P$ be given. Then $t \leq ts$ which implies $M(ts) \leq M(t)$. Thus there exists a $t_{i}$ such that $M(ts)M(t_{i})\neq 0$. This implies that $ts$ and $t_{i}$ have a common upper bound in $P$. In other words, $(ts)^{-1}t_{i}=s^{-1}t^{-1}t_{i} \in PP^{-1}$. Thus $s$ and $t^{-1}t_{i}$ have a common upper bound in $P$. This proves the claim.

By assumption it follows that $\prod_{i=1}^{n}(1-L(t^{-1}t_{i}))=0$ where $L(s):=V(s)V(s)^{*}$. Now multiplying this equality on the left by $V(t)$ and on the right by $V(t)^{*}$, we get 
\begin{align*}
  \prod_{i=1}^{n}(V(t)V(t)^{*}-V(t)V(t^{-1}t_{i})V(t^{-1}t_{i})^{*}V(t)^{*}) &=0\\
  \prod_{i=1}^{n}(V(t)V(t)^{*}-V(t_{i})V(t_{i})^{*})&=0
\end{align*}
 But this is equivalent to $\displaystyle \sup_{i} L(t_{i})=L(t)$. This completes the proof. \hfill $\Box$
\begin{rmrk}
 The relations $\prod_{x \in F}(1-V(t)V(t)^{*})=0$ for $F \in \mathcal{F}$ are  the boundary relations considered in \cite{Crisp}.
\end{rmrk}

Let $Q_{\mathbb{N}}$ be the $C^{*}$-subalgebra of $U[\mathbb{Z}]$ generated by $u$ and $\{s_{m}:m>0\}$. In \cite{Cuntz}, it was proved that $Q_{\mathbb{N}}$ is simple and purely infinite. Moreover in \cite{Cuntz}, it was shown that $U[\mathbb{Z}]$ is isomorphic to a crossed product of $Q_{\mathbb{N}}$ with $\mathbb{Z}/2\mathbb{Z}$. Let 
\begin{equation*}
 P_{\mathbb{N}}:=\Big\{ \begin{bmatrix}
                     1 & 0 \\
                     k & m \\
                    \end{bmatrix}: k\in \mathbb{N} \text{~and~} m \in \mathbb{N}^{\times} \Big \}
\end{equation*}
Note that $P_{\mathbb{N}}$ is a semigroup of $P_{\mathbb{Q}}$.
\begin{rmrk}
In \cite{Laca1}, it was proved that $(P_{\mathbb{Q}},P_{\mathbb{N}})$ is a quasi-lattice ordered group. Moreover it was shown in \cite{Laca1} that Nica-covariance together with boundary relations is equivalent to Cuntz-Li relations and the universal $C^{*}$-algebra made out of Nica-covariant representations satisfying the boundary relations is in fact $Q_{\mathbb{N}}$.
\end{rmrk}

\section{The Cuntz-Li algebra for a general integral domain}
  We end this article by giving a few remarks of how to adapt the  analysis in Section $1-6$ for a general integral domain $R$. Now Let $R$ be an integral domain such that $R/mR$ is finite for every non-zero $m \in R$. We also assume that $R$ is countable and $R$ is not a field.
\begin{dfn}[\cite{Cuntz-Li}]
  Let $U[R]$ be the universal $C^*$-algebra generated by a set of
unitaries $\{u^{n}:n \in R\}$ and a set of isometries $\{s_{m}: m \in
R^{\times} \}$ satisfying the following relations.
\begin{align*}
 s_{m}s_{n} &= s_{mn} \\
 u^{n}u^{m}&=u^{n+m} \\
 s_{m}u^{n}&=u^{mn}s_{m} \\
\sum_{n \in R/mR}u^{n}e_{m}u^{-n}&=1
 \end{align*}
where $e_{m}$ denotes the final projection of $s_{m}$.
\end{dfn}
 Now the problem is the product $u^{r}e_{m}u^{-r}u^{s}e_{n}u^{-s}$ may not be of the form $u^{k}e_{c}u^{-k}$ for some $k$ and $c$. Nevertheless it will be in the linear span of $\{u^{k}e_{mn}u^{-k}:k \in R/(mn)\}$. 
Let $P$ denote the set of projections in $U[R]$ which is in the linear span of $\{u^{r}e_{m}u^{-r}: r \in R/(m)\}$ for some $m$. Explicity, a projection $e \in U[R]$ is in $P$ if and only if there exists an $m \in R^{\times}$ and $a_{r} \in \{0,1\}$ such that $f=\sum_{r}a_{r}u^{r}e_{m}u^{-r}$.

Now it is easy to show that  $P$ is a commutative semigroup of projections containing $0$. Moreover $P$ is invariant under conjugation by $u^{r}$, $s_{m}$ and $s_{m}^{*}$. One can prove the following Proposition just as in the case when $R=\mathbb{Z}$.

\begin{ppsn}
\label{Inverse semigroup1}
 Let $T:=\{s_{m}^{*}u^{n}eu^{n^{'}}s_{m^{'}}:e \in P, m,m^{'} \neq 0, n ,n^{'} \in
R\}$. Then $T$ is an inverse semigroup of partial
isometries. Moreover the set of projections in $T$ coincide with $P$. Also the linear span of $T$ is dense in $U[R]$.
\end{ppsn}
 
Let $\hat{R}:=\{(r_{m})\in \prod R/(m): r_{mk}=r_{m} \textrm{~in~}R/(m)\}$ be the profinite completion of the ring $R$. For $r \in \hat{R}$, define 
\begin{displaymath}
 A_{r}:=\{f \in P: f \geq u^{r_{m}}e_{m}u^{-r_{m}} \textrm{~for some ~} m\}
\end{displaymath}
Then $A_{r}$ is an ultrafilter for every $r \in \hat{R}$ and  the map $r \to A_{r}$ is a topological isomorphism from  $\hat{R}$ to $\hat{P}_{tight}$.

 Let $Q(R)$ be the field of fractions of $R$. For $m\neq 0$, let $\mathcal{R}_{m}:=\hat{R}$. For every $\ell \neq 0$, let $\phi_{m\ell,m}:\mathcal{R}_{m} \to \mathcal{R}_{\ell m}$ be the map defined by mulitplication by $\ell$. Then $\phi_{m\ell,m}$ is only an additive homomorphism and it does not preserve the multiplication. We let $\mathcal{R}$ be the inductive limit of $(\mathcal{R}_{m},\phi_{m\ell,m})$. Then $\mathcal{R}$ is an abelian group and $\hat{R}$ is a subgroup of $\mathcal{R}$ via the inclusion $\mathcal{R}_{1}\subset \mathcal{R}$. Note that $\mathcal{R}$ is a locally compact Hausdorff space. Moreover  the group 
\[                                                                                                                                                                                                                                                                                       
 P_{Q(R)}:=\big\{\begin{bmatrix}
                          1 & 0 \\
                          b & a                                                                                                                                                                                                                                                                                                                                                                                                                                                                                                                                                                                                                                                                                                                                                                                   
                                                                                                                                                                                                                                                                                                                                                                                                                                                                                                                                                                                                                                                                                                                                                                                             \end{bmatrix}: a \in Q(R)^{\times},b \in Q(R) \big \}                                                                                                                                                                                                                                                                                                                                                                                                                                                                                                                                                                                                                               \]
 acts on $\mathcal{R}$ by affine transformations. The action is descibed explicitly by the following formula. For $x \in \mathcal{R}_{p}$
\begin{displaymath}
 \begin{bmatrix}
  1 & 0 \\
  \frac{n}{m^{'}} & \frac{m}{m^{'}}
 \end{bmatrix}x = mx+np \in \mathcal{R}_{m^{'}p}
\end{displaymath}

One can check that the above formula defines an action of $P_{Q(R)}$ on $\mathcal{R}$. Let $\mathcal{G}_{tight}$ be the tight groupoid associated to the inverse semigroup $T$ defined in Proposition \ref{Inverse semigroup1}. Then as in the case when $R=\mathbb{Z}$, we have the following theorem.

\begin{thm}
\label{main theorem}
 Let $\phi:\mathcal{R}\times P_{Q(R)}|_{\hat{R}} \to \mathcal{G}_{tight}$ be the map defined by 
\begin{displaymath}
\phi\Big( \big(r,\begin{bmatrix}
                                                                           1 & 0 \\
                                                                         \frac{k}{m} & \frac{n}{m}
                                                                          \end{bmatrix}
\big)\Big)=
[(r,s_{m}^{*}u^{k}s_{n})]
\end{displaymath}
 Then $\phi$ is a topological groupoid isomorphism. Moreover the $C^{*}$-algebra $U[R]$ is isomorphic to the full ( and the reduced) $C^{*}$-algebra of the groupoid $\mathcal{R}\times P_{Q(R)}|_{\hat{R}}$.
 \end{thm}
We end this article by showing that $U[R]$ is simple.
\begin{ppsn}[\cite{Cuntz-Li}]
 The $C^{*}$-algebra $U[R]$ is simple.
\end{ppsn}
\textit{Proof.} Let us denote the groupoid $\mathcal{R}\times P_{Q(R)}|_{\hat{R}}$ by $\mathcal{G}$. As in Proposition \ref{simplicity}, we need to show that $\mathcal{G}$ is minimal and topologically principal. The proof of the minimality of $\mathcal{G}$ is exactly similar to that in Proposition \ref{simplicity}. We now show that 
$\mathcal{G}$ is topologically principal. For $g \in P_{Q(R)} \backslash \{1\}$, let us denote the set of fixed points of $g$ in $\hat{R}$ by $F_{g}$. It follows from Baire category theorem that $\mathcal{G}$ is topologically principal if and only if $F_{g}$ has empty interior for every $g\neq 1$.

Let $g=\begin{bmatrix}
        1 & 0 \\
        \frac{k}{m} & \frac{n}{m}
       \end{bmatrix}$ be a non-identity element in $P_{Q(R)}$. Suppose that $F_{g}$ contains a non-empty open set say $U$. Now note that $R$ is dense in $\hat{R}$. Thus $U \cap R$ is non-empty. Moreover 
       $U\cap R$ is  infinite. Let $r_{1},r_{2}$ be two distinct points of $R$ in $U$. Since $r_{1},r_{2} \in F_{g}$, it follows that $mr_{1}-k=nr_{1}$ and $mr_{2}-k=nr_{2}$. Thus we have $(m-n)r_{1}=k=(m-n)r_{2}$. This forces $m=n$ and $k=0$. This is a contradiction to the fact that $g\neq 1$. Thus for every $g\neq 1$, $F_{g}$ has empty interior which in turn implies that $\mathcal{G}$ is topologically principal. This completes the proof. \hfill $\Box$

\begin{rmrk}
 In \cite{Quigg-Landstad}, Cuntz-Li type relations arising out of a semidirect product $N \rtimes H$ where $N$ is a normal subgroup and $H$ is an abelian group satisfying certain hypothesis were considered. It was shown in \cite{Quigg-Landstad} that the universal $C^{*}$-algebra generated by the Cuntz-Li type relations is isomorphic to a corner of a crossed product algebra. It is possible to apply inverse semigroups and tight representations to reconstruct this result. The details will be spelt out elsewhere.
\end{rmrk}

\bibliography{references}
\bibliographystyle{amsalpha}

\nocite{Cuntz-Li-1}
\nocite{Cuntz}
\nocite{Muhly}
\nocite{Crisp}
\nocite{Laca1}
\nocite{Nica1}

\end{document}